\documentclass[11pt]{article}
\usepackage{mathrsfs}
\usepackage{amssymb}
\usepackage{amsmath}
\usepackage[all]{xy}

\def\Tr{\mathop{\rm Tr}\nolimits}

\title{\Large \bf Homological Behavior of Auslander's $k$-Gorenstein Rings
}
\author{Zhaoyong Huang\thanks{{\it E-mail address}: huangzy@nju.edu.cn}
{\small \ and}
Hourong Qin\thanks{{\it E-mail address}: hrqin@nju.edu.cn}\\
{\footnotesize \it Department of Mathematics, Nanjing University, Nanjing 210093, Jiangsu Province, P. R. China}\\ \\}
\usepackage{amssymb}
\date{}
\begin{document}
\baselineskip=18pt \maketitle

\begin{abstract} In this paper we mainly study the homological properties of dual
modules over $k$-Gorenstein rings. For a right quasi $k$-Gorenstein
ring $\Lambda$, we show that the right self-injective dimension of
$\Lambda$ is at most $k$ if and only if each $M \in$mod $\Lambda$
satisfying the condition that Ext$_{\Lambda}^i(M, \Lambda)=0$ for
any $1\leq i \leq k$ is reflexive. For an $\infty$-Gorenstein ring,
we show that the big and small finitistic dimensions and the
self-injective dimension of $\Lambda$  are identical. In addition,
we show that if $\Lambda$ is a left quasi $\infty$-Gorenstein ring
and $M\in$mod $\Lambda$ with grade$M$ finite, then
Ext$_{\Lambda}^i($Ext$_{\Lambda ^{op}}^i($Ext$_{\Lambda}^{{\rm
grade}M}(M, \Lambda), \Lambda), \Lambda)=0$ if and only if
$i\neq$grade$M$. For a 2-Gorenstein ring $\Lambda$, we show that a
non-zero proper left ideal $I$ of $\Lambda$ is reflexive if and only
if $\Lambda /I$ has no non-zero pseudo-null submodule.
\end{abstract}

\vspace{0.5cm}

\centerline{\large \bf 1. Introduction}

\vspace{0.2cm}

Throughout this paper, $\Lambda$ is a left and right Noetherian
ring. We use mod $\Lambda$ to denote the category of finitely
generated left $\Lambda$-modules. Let $M\in$mod $\Lambda$. We use
pd$_{\Lambda}M$, fd$_{\Lambda}M$ and id$_{\Lambda}M$ to denote the
projective, flat and injective dimensions of $M$, respectively. For
a non-negative integer $i$, recall from [AR4] that the {\it grade}
of $M$, denoted by grade$M$, is said to be at least $i$ if
Ext$_{\Lambda}^{j}(M, \Lambda)=0$ for any $0\leq j<i$; and the {\it
strong grade} of $M$, denoted by s.grade$M$, is said to be at least
$i$ if grade$X\geq i$ for each submodule $X$ of $M$. In addition, we
fix minimal injective resolutions
$$0\to {_{\Lambda}\Lambda} \to I^{'}_0 \to I^{'}_1 \to \cdots \to
I^{'}_i \to \cdots$$ and $$0\to {\Lambda_{\Lambda}} \to I_0 \to I_1
\to \cdots \to I_i \to \cdots$$ of ${_{\Lambda}\Lambda}$ and
${\Lambda_{\Lambda}}$, respectively.

Recall from [FGR] that $\Lambda$ is called a $k$-{\it Gorenstein
ring} if fd$_{\Lambda}I^{'}_i\leq i$ for any $0\leq i \leq k-1$, and
$\Lambda$ is called an {\it $\infty$-Gorenstein ring} if it is
$k$-Gorenstein for all $k$. Auslander gave some equivalent
characterizations of $k$-Gorenstein rings in terms of the strong
grade and the flat dimension of injective modules as follows.

\vspace{0.2cm}

{\bf Auslander's Theorem} ([FGR, Theorem 3.7]) {\it The following
statements are equivalent.}

(1) $\Lambda$ {\it is a} $k$-{\it Gorenstein ring}.

(2) fd$_{\Lambda ^{op}}I_i\leq i$ {\it for any} $0\leq i \leq k-1$.

(3) s.gradeExt$_{\Lambda}^i(M, \Lambda)\geq i$ {\it for any} $M
\in$mod $\Lambda$ {\it and} $1\leq i \leq k$.

(4) s.gradeExt$_{\Lambda ^{op}}^i(N, \Lambda)\geq i$ {\it for any}
$N \in$mod $\Lambda ^{op}$ {\it and} $1\leq i \leq k$.

\vspace{0.2cm}

The Auslander's theorem shows that the notion of $k$-Gorenstein rings is
left-right symmetric. The properties of $k$-Gorenstein rings and
related rings have been studied by many authors (see [AB], [AR3],
[AR4], [Bj], [BjE], [CSS], [FGR], [HN], [Hu1], [Hu3], [Hu4], [I2],
[IS], [Iy], [M], and so on). In this paper we mainly investigate the
homological properties of dual modules over $k$-Gorenstein rings and
related rings.

It was showed in [HuT, Theorem 2.2] that if id$_{\Lambda
^{op}}\Lambda\leq k$, then each $M \in$mod $\Lambda$ satisfying the
condition that Ext$_{\Lambda}^i(M, \Lambda)=0$ for any $1\leq i \leq
k$ is reflexive. In Section 2, we show that the converse of this
result holds true for right quasi $k$-Gorenstein rings (see Section 2 for
the definition of such rings).
It was showed in [Zi] that the big and small
finitistic dimensions of $\Lambda$ are usually different even when
$\Lambda$ is an Artinian algebra. The other aim of Section 2 is to
show that for a $(k+1)$-Gorenstein ring $\Lambda$, the small
finitistic dimension of $\Lambda$ is at most $k$ if and only if
id$_{\Lambda}\Lambda\leq k$, which yields that for an
$\infty$-Gorenstein ring $\Lambda$, the big and small finitistic
dimensions of $\Lambda$ and id$_{\Lambda}\Lambda$ are identical. As
a consequence, we get some equivalent versions of the Nakayama
conjecture.

In Section 3, we get some properties of the grade of modules over quasi
$\infty$-Gorenstein rings. We prove that if $\Lambda$ is a left
quasi $\infty$-Gorenstein ring and $M\in$mod $\Lambda$ with grade$M$
finite, then Ext$_{\Lambda}^i($Ext$_{\Lambda
^{op}}^i($Ext$_{\Lambda}^{{\rm grade}M}(M, \Lambda), \Lambda),
\Lambda)=0$ if and only if $i\neq$grade$M$. As an application of
this result, we show that if $\Lambda$ is a left and right quasi
Auslander-Gorenstein ring with id$_{\Lambda}\Lambda$=id$_{\Lambda
^{op}}\Lambda =t$, then the functors Ext$_{\Lambda}^t(\ , \Lambda)$
and Ext$_{\Lambda ^{op}}^t(\ , \Lambda)$ give a duality between $\{
M\in {\rm mod}\ \Lambda | {\rm grade}M=t\}$ and $\{ N\in {\rm mod}\
\Lambda ^{op}| {\rm grade}N=t\}$. This generalizes a result of
Iwanaga in [I2].

In Section 4, we give a criterion for judging when a finitely
generated torsionless module is reflexive. As a consequence
of this criterion, we prove that for a 2-Gorenstein ring
$\Lambda$, a non-zero proper left ideal $I$ of $\Lambda$ is
reflexive if and only if $\Lambda /I$ has no non-zero pseudo-null
submodule. This generalizes a result of Coates, Schneider and
Sujatha in [CSS]. We also study the following question:
For a positive integer $k$, when is each $k$-torsionfree module
in mod $\Lambda$ projective? We prove
that the answer to this question is affirmative if gl.dim$\Lambda$
(the global dimension of $\Lambda$) is at most $k$.

\vspace{0.5cm}

\centerline{\large \bf 2. Reflexive modules and homological
dimensions}

\vspace{0.2cm}

{\bf Lemma 2.1} {\it The following statements are equivalent for a
positive integer} $k$.

(1)  fd$_{\Lambda}I_i^{'}\leq i+1$ {\it for any} $0\leq i \leq k-1$.

(2) s.gradeExt$_{\Lambda ^{op}}^{i+1}(N, \Lambda)\geq i$ {\it for
any} $N\in$ mod $\Lambda ^{op}$ {\it and} $1\leq i\leq k$.

(3) gradeExt$_{\Lambda}^i(M, \Lambda)\geq i$ {\it for any} $M\in$
mod $\Lambda$ {\it and} $1\leq i\leq k$.

\vspace{0.2cm}

{\it Proof.} By [AR4, Theorem 0.1] and [HN, Theorem 4.1].
\hfill{$\square$}

\vspace{0.2cm}

{\bf Definition 2.2}$^{{\rm [Hu3, Hu4]}}$ $\Lambda$ is called a {\it
left quasi} $k$-{\it Gorenstein ring} if one of the equivalent
conditions in Lemma 2.1 is satisfied. $\Lambda$ is called a {\it
left quasi $\infty$-Gorenstein ring} if it is left quasi
$k$-Gorenstein for all $k$. Dually, the notions of
{\it right quasi} $k$-{\it Gorenstein rings} and {\it right quasi
$\infty$-Gorenstein rings} are defined.

\vspace{0.2cm}

{\it Remark.} A $k$-Gorenstein ring clearly is a left and right
quasi $k$-Gorenstein ring. However, the interesting property of the
left-right symmetry enjoyed by $k$-Gorenstein rings fails for quasi
$k$-Gorenstein rings (see [Hu1]).

\vspace{0.2cm}

Let $M\in$mod $\Lambda$ and $k$ a positive integer. Recall from [H1]
that the {\it reduced grade} of $M$, denoted by r.grade$M$, is said to be
at least $k$ if Ext$_{\Lambda}^{j}(M, \Lambda)=0$ for any $1\leq
j<k$. It is trivial that grade$M \geq k$ implies r.grade$M \geq k$.
Let $\sigma _{M}: M \to M^{**}$ via $\sigma _{M}(x)(f)=f(x)$ for any
$x\in M$ and $f\in M^*$ be the canonical evaluation homomorphism,
where $(\ )^*=$Hom$(\ , \Lambda)$. Recall that $M$ is called {\it
torsionless} if $\sigma _{M}$ is a monomorphism; and $M$ is called
{\it reflexive} if $\sigma _{M}$ is an isomorphism. It was showed in
[HuT, Theorem 2.2] that if id$_{\Lambda ^{op}}\Lambda \leq k$, then
each module in mod $\Lambda$ with reduced grade at least $k+1$ is
reflexive. However, we do not know whether the converse holds true.
Jans showed in [J2, Theorem 5.1] that the converse holds true when
$k=1$. One of the main results in this section is the following

\vspace{0.2cm}

{\bf Theorem 2.3} {\it Let} $\Lambda$ {\it be a right quasi}
$k$-{\it Gorenstein ring. Then the following statements are
equivalent.}

(1) id$_{\Lambda ^{op}}\Lambda \leq k$.

(2) {\it Each module in} mod $\Lambda$ {\it with reduced grade at
least} $k+1$ {\it is reflexive}.

(3) {\it Each module in} mod $\Lambda$ {\it with reduced grade at
least} $k+1$ {\it is torsionless}.

\vspace{0.2cm}

Before proving this theorem, we recall some notions from [AR2] and
[J2].\linebreak
Let $M\in$mod $\Lambda$ and $\mathcal{D}$ be a full
subcategory of mod $\Lambda$. A homomorphism $M\buildrel {f} \over
\to D$ in mod $\Lambda$ with $D\in \mathcal{D}$ is called a left
$\mathcal{D}$-{\it approximation} of $M$ if Hom$_{\Lambda}(f, D')$
is epic for any $D'\in \mathcal{D}$ (see [AR2]). Let $X\in$mod
$\Lambda$ and $Y\in$mod $\Lambda ^{op}$. A monomorphism $X^{**}
\buildrel {\rho ^*}\over \to Y^*$ is called a {\it double dual
embedding} if it is the dual of an epimorphism $Y \buildrel
{\rho}\over \to X^*$. For a positive integer $k$, a torsionless
module $T_k\in$mod $\Lambda$ is said to be of $D$-{\it class} $k$ if
it can be fitted into a sequence of $k-1$ exact sequences of the
form:
$$\xymatrix{& & & & 0\ar[r] & T_{k-1}^{**}\ar[r]
& P_{k-1}\ar[r] & T_k \ar[r] & 0 \\
& & & \cdots \ar[r] & P_{k-2}\ar[r] & T_{k-1}\ar[r]\ar[u]_{\sigma
_{T_{k-1}}} & 0 & & \\
& & & \cdots & & & & & \\
& & 0\ar[r] & T_{2}^{**}\ar[r] & \cdots & &
& & \\
0 \ar[r] & T_1^{**} \ar[r] & P_1 \ar[r] & T_2 \ar[r]\ar[u]_{\sigma
_{T_{2}}} & 0 & & & & }
$$ where each $P_i\in$mod $\Lambda$ is projective
and the horizontal monomorphisms are double dual embeddings. Any
torsionless module in mod $\Lambda$ is said to be of $D$-{\it class
1} (see [J2]). It follows from [J2, p.335] that a torsionless module
in mod $\Lambda$ is of $D$-class $k$ if its reduced grade is at
least $k$. The following two results are cited from [J2], which play
a crucial role in proving Theorem 2.3.

\vspace{0.2cm}

{\bf Lemma 2.4} ([J2, Theorem 4.3]) id$_{\Lambda ^{op}}\Lambda \leq
k$ {\it if and only if each module of} $D$-{\it class} $k$ {\it in}
mod $\Lambda$ {\it is reflexive}.

\vspace{0.2cm}

The proof of [J2, Theorem 3.1] can be applied to obtain the
following more general result.

\vspace{0.2cm}

{\bf Lemma 2.5} {\it If} $\Lambda$ {\it is a right quasi} $k$-{\it
Gorenstein ring, then the reduced grade of each module of} $D$-{\it
class} $k$ {\it in} mod $\Lambda$ {\it is at least} $k$.

\vspace{0.2cm}

{\it Proof of Theorem 2.3.} $(1)\Rightarrow (2)$ follows from [HuT,
Theorem 2.2], and $(2)\Rightarrow (3)$ is trivial.

$(3)\Rightarrow (1)$ Assume that $M \in$mod $\Lambda$ is of
$D$-class $k$. Then $M$ is torsionless and $\sigma _M$ is a
monomorphism. By Lemma 2.5, we have that r.grade$M\geq k$.

Let $P\buildrel {f} \over \to M^* \to 0$ be exact in mod $\Lambda
^{op}$ with $P$ projective. Put $g=f^* \sigma _M$ and $X=$Coker$g$.
It is not difficult to verify that the monomorphism $g: M\to P^*$ is
a left $\mathcal{P}^0(\Lambda)$-approximation of $M$, where
$\mathcal{P}^0(\Lambda)$ denotes the full subcategory of mod
$\Lambda$ consisting of projective modules. Then r.grade$X\geq 2$
and so r.grade$X\geq k+1$. Thus by (3) we have that $X$ is torsionless
and $\sigma _X$ is a monomorphism. Furthermore, we have the
following commutative diagram with exact rows:
$$\xymatrix{0\ar[r] & M\ar[r]^g \ar[d]^{\sigma _M}
& P^* \ar[r] \ar[d]^{\sigma _{P^*}}_{\cong} & X \ar[r]
\ar[d]^{\sigma _X} & 0\\
0\ar[r] & M^{**}\ar[r]^{g^{**}} & P^{***} \ar[r] & X^{**} & }
$$ Then by the snake lemma we have that $\sigma _M$ is an isomorphism and $M$ is
reflexive. So id$_{\Lambda ^{op}}\Lambda
\leq k$ by Lemma 2.4. \hfill{$\square$}

\vspace{0.2cm}

The following lemma improves [Hu2, Corollary 3].

\vspace{0.2cm}

{\bf Proposition 2.6} {\it Let} $\Lambda$ {\it be a left and right
Artinian ring. If} id$_{\Lambda ^{op}}\Lambda =k$ {\it and}
fd$_{\Lambda ^{op}}\bigoplus _{i=0}^{k-2}I_i<\infty$, {\it then}
id$_{\Lambda}\Lambda =k$.

\vspace{0.2cm}

{\it Proof.} Assume that fd$_{\Lambda ^{op}}\bigoplus
_{i=0}^{k-2}I_i=r(<\infty)$. By [Iy, 6.1(1)], we have that
\linebreak
s.gradeExt$_{\Lambda}^{r+1}(M, \Lambda)\geq k-1$ for any $M \in$mod
$\Lambda$. Then id$_{\Lambda}\Lambda \leq (r+1)+k-1=r+k$ by [Hu2,
Theorem]. It follows from [Z, Lemma A] that id$_{\Lambda}\Lambda
=$id$_{\Lambda ^{op}}\Lambda =k$. \hfill{$\square$}

\vspace{0.2cm}

By Proposition 2.6 and the left-right symmetry of $k$-Gorenstein
rings, we immediately have the following

\vspace{0.2cm}

{\bf Corollary 2.7} {\it Let} $\Lambda$ {\it be a left and right
Artinian ring. If} $\Lambda$ {\it is} $(k-1)$-{\it Gorenstein, then}
id$_{\Lambda ^{op}}\Lambda \leq k$ {\it if and only if}
id$_{\Lambda}\Lambda \leq k$.

\vspace{0.2cm}

The following result generalizes [AR3, Corollary 5.5(b)].

\vspace{0.2cm}

{\bf Corollary 2.8} {\it Let} $\Lambda$ {\it be a left and right
Artinian ring. If} $\Lambda$ {\it is} $\infty$-{\it Gorenstein,
then} id$_{\Lambda ^{op}}\Lambda =$id$_{\Lambda}\Lambda$.

\vspace{0.2cm}

{\bf Corollary 2.9} {\it Let} $\Lambda$ {\it be a left and right
Artinian ring. If} $\Lambda$ {\it is} $k$-{\it Gorenstein, then the
following statements are equivalent}.

(1) id$_{\Lambda ^{op}}\Lambda \leq k$.

(2) {\it Each module in} mod $\Lambda$ {\it with reduced grade at
least} $k+1$ {\it is reflexive}.

(3) {\it Each module in} mod $\Lambda$ {\it with reduced grade at
least} $k+1$ {\it is torsionless}.

(1)$^{op}$ id$_{\Lambda}\Lambda \leq k$.

(2)$^{op}$ {\it Each module in} mod $\Lambda ^{op}$ {\it with
reduced grade at least} $k+1$ {\it is reflexive}.

(3)$^{op}$ {\it Each module in} mod $\Lambda ^{op}$ {\it with
reduced grade at least} $k+1$ {\it is torsionless}.

\vspace{0.2cm}

{\it Proof.} By Theorem 2.3 and Corollary 2.7. \hfill{$\square$}

\vspace{0.2cm}

Recall that the {\it big finitistic dimension} of $\Lambda$, denoted
by Fin.dim$\Lambda$, is defined to be sup$\{$pd$_{\Lambda}M|M$ is a
left $\Lambda$-module with pd$_{\Lambda}M<\infty \}$; and the {\it
small finitistic dimension} of $\Lambda$, denoted by
fin.dim$\Lambda$, is defined to be sup$\{$pd$_{\Lambda}M|M\in$mod
$\Lambda$ with pd$_{\Lambda}M<\infty \}$. According to [Zi], the big
and small finitistic dimensions are not identical in general even
for Artinian algebras. The other aim of this section is to show that
for an $\infty$-Gorenstein ring $\Lambda$, these two dimensions and
id$_{\Lambda}\Lambda$ are identical. As an application of this
result, we get some equivalent versions of the Nakayama conjecture.

We begin with the following easy observation.

\vspace{0.2cm}

{\bf Lemma 2.10} {\it Let} $M \in$mod $\Lambda$ {\it with} r.grade$M
\geq k+1$. {\it If} pd$_{\Lambda}M\leq k$, {\it then} $M$ {\it is
projective.}

\vspace{0.2cm}

{\it Proof.} Consider the projective resolution of $M$ in mod
$\Lambda$. Then we get easily the assertion. \hfill{$\square$}

\vspace{0.2cm}

{\bf Lemma 2.11} ([B2, Proposition 4.3]) fin.dim$\Lambda
\leq$Fin.dim$\Lambda \leq$id$_{\Lambda}\Lambda$.

\vspace{0.2cm}

Let $M\in$mod $\Lambda$ and $k$ a positive integer. Recall that $M$
is called a {\it $k$-syzygy module} if there exists an exact
sequence $0 \to M \to Q_{0} \to Q_{1} \to \cdots \to Q_{k-1}$ in mod
$\Lambda$ with all $Q_{i}$ projective. On the other hand, assume
that
$$P_{1} \buildrel {f} \over \longrightarrow P_{0} \to M \to 0$$ is a
projective presentation of $M$ in mod $\Lambda$. Then we get an
exact sequence: $$0 \to M^* \to P_{0}^* \buildrel {f^*} \over
\longrightarrow P_{1}^* \to \Tr M \to 0$$ in mod $\Lambda^{op}$, where $\Tr M={\rm
Coker}f^*$ is the {\it transpose} of $M$. $M$ is called a {\it
$k$-torsionfree module} if r.grade$\Tr M\geq k+1$; and $M$ is called
an $\infty$-{\it torsionfree module} if r.grade$\Tr M\geq k$ for all
$k$ (see [AB]). We remark that the transpose of $M$ depends on the
choice of the projective resolution of $M$, but it is unique up to
projective equivalence. So the notions of $k$-torsionfree modules
and $\infty$-torsionfree modules are well defined.

\vspace{0.2cm}

{\bf Lemma 2.12} ([AB, Proposition 2.6]) {\it Let} $M \in$mod
$\Lambda$. {\it Then we have the following exact sequences:}
$$0\to {\rm Ext}_{\Lambda ^{op}}^1(\Tr M, \Lambda)\to M \buildrel {\sigma _M}
\over \longrightarrow M^{**} \to {\rm Ext}_{\Lambda ^{op}}^2(\Tr M,
\Lambda) \to 0,$$
$$0\to {\rm Ext}_{\Lambda}^1(M, \Lambda)\to \Tr M \buildrel {\sigma _{\Tr M}}
\over \longrightarrow (\Tr M)^{**} \to {\rm Ext}_{\Lambda}^2(M,
\Lambda) \to 0.$$

\vspace{0.2cm}

It follows from Lemma 2.12 that a module in mod $\Lambda$ is
torsionless (resp. reflexive) if and only if it is 1-torsionfree
(resp. 2-torsionfree). The following observation is well known,
which is an immediate consequence of Lemma 2.12.

\vspace{0.2cm}

{\bf Proposition 2.13} {\it A module} $M$ {\it in} mod $\Lambda$
{\it is projective if and only if} $\Tr M$ {\it is projective in}
mod $\Lambda ^{op}$.

\vspace{0.2cm}

{\it Proof.} The necessity is trivial, so it suffices to prove the
sufficiency. Let $M \in$mod $\Lambda$ with $\Tr M$
projective. Then $M^*$ is projective. So by Lemma 2.12 we have
that $\sigma _M$ is an isomorphism and $M(\cong M^{**})$ is
projective. \hfill{$\square$}

\vspace{0.2cm}

{\bf Lemma 2.14} {\it Consider the following conditions.}

(1) {\it Each module in} mod $\Lambda ^{op}$ {\it is reflexive}.

(2) {\it Each module in} mod $\Lambda ^{op}$ {\it is torsionless}.

(3) $\Lambda$ {\it is self-injective}.

(4) fin.dim$\Lambda =0$.

(5) {\it A module} $N \in$mod $\Lambda ^{op}$ {\it is projective
provided} $N^*$ {\it is projective}.

{\it We have} $(1)\Leftrightarrow (2)\Leftrightarrow (3)\Rightarrow
(4)\Leftrightarrow (5)$. {\it If} $\Lambda$ {\it is a 1-Gorenstein
ring, then all of these conditions are equivalent.}

\vspace{0.2cm}

{\it Proof.} $(1) \Rightarrow (2)$ is trivial. $(1)\Leftrightarrow
(3)$ and $(3)\Rightarrow (4)$ follow from [J1, Corollary 1.2] and
Lemma 2.11, respectively.

$(2) \Rightarrow (3)$ Let $M\in$mod $\Lambda$. Then $\Tr M\in$mod
$\Lambda ^{op}$ is torsionless by (2). It follows from Lemma 2.12
that ${\rm Ext}_{\Lambda}^1(M, \Lambda)\cong \sigma _{\Tr M}=0$ and
$\Lambda$ is self-injective.

$(4)\Rightarrow (5)$ Let $N\in$mod $\Lambda ^{op}$ with $N^*$
projective. Then $\Tr N \in$mod $\Lambda$ and pd$_{\Lambda}\Tr N\leq
2$. By (4), $\Tr N$ is projective. It follows from Proposition 2.13
that $N$ is projective.

$(5)\Rightarrow (4)$ Let $N\in$mod $\Lambda ^{op}$ with $N^*=0$,
Then $N$ is projective by (5) and so $N=0$. It follows from [B1,
Corollary 5.6 and Theorem 5.4] that fin.dim$\Lambda =0$.

Now let $\Lambda$ be a 1-Gorenstein ring. We will prove $(4)
\Rightarrow (3)$. Suppose fin.dim$\Lambda =0$. Because $\Lambda$ is
1-Gorenstein, for any $M\in$mod $\Lambda$ we have that
s.gradeExt$_{\Lambda}^1(M, \Lambda)\geq 1$ and $[{\rm
Ext}_{\Lambda}^1(M, \Lambda)]^*=0$. Then by [B1, Corollary 5.6 and
Theorem 5.4], ${\rm Ext}_{\Lambda}^1(M, \Lambda)=0$. It implies that
$\Lambda$ is self-injective. \hfill{$\square$}

\vspace{0.2cm}

{\bf Theorem 2.15} {\it Let} $k$ {\it be a non-negative integer and}
$\Lambda$ {\it a} $(k+1)$-{\it Gorenstein ring with}
fin.dim$\Lambda=k$. {\it Then} id$_{\Lambda}\Lambda=k$.

\vspace{0.2cm}

{\it Proof.} The case for $k=0$ follows from Lemma 2.14. Now suppose
that $k\geq 1$ and $M$ is any module in mod $\Lambda$. Since
$\Lambda$ is $(k+1)$-Gorenstein, s.gradeExt$_{\Lambda}^{k+1}(M,
\Lambda)\geq k+1$. Let
$$Q_{k+1} \to Q_k \to \cdots \to Q_1 \to Q_0 \to
{\rm Ext}_{\Lambda}^{k+1}(M, \Lambda) \to 0$$ be a projective
resolution of Ext$_{\Lambda}^{k+1}(M, \Lambda)$ in mod $\Lambda
^{op}$. Put $N_k=$Coker$(Q_{k+1} \to Q_k)$. Notice that a $k$-syzygy
module is $k$-torsionfree by [AR4, Proposition 1.6], so $N_k$ is
$k$-torsionfree and r.grade$\Tr N_k\geq k+1$. On the other hand, we
have the following exact sequence:
$$0 \to Q_0^* \to Q_1^* \to \cdots \to Q_k^* \to Q_{k+1}^* \to \Tr N_k
\to 0$$ in mod $\Lambda$. Then pd$_{\Lambda}\Tr N_k \leq k+1$. But
fin.dim$\Lambda=k$, so pd$_{\Lambda}\Tr N_k \leq k$ and hence $\Tr
N_k$ is projective by Lemma 2.10. It follows from Proposition 2.13
that $N_k$ is projective and thus pd$_{\Lambda
^{op}}$Ext$_{\Lambda}^{k+1}(M, \Lambda)\leq k$. Now assume that
pd$_{\Lambda ^{op}}$Ext$_{\Lambda}^{k+1}(M, \Lambda)=t(\leq k)$. If
Ext$_{\Lambda}^{k+1}(M, \Lambda)\neq 0$, then it is not difficult to
verify that ${\rm Ext}_{\Lambda ^{op}}^t({\rm
Ext}_{\Lambda}^{k+1}(M, \Lambda), \Lambda)\neq 0$, which contradicts
with that $\Lambda$ is $(k+1)$-Gorenstein. Thus
Ext$_{\Lambda}^{k+1}(M, \Lambda)=0$ and therefore
id$_{\Lambda}\Lambda \leq k$. It follows from Lemma 2.11 that
id$_{\Lambda}\Lambda =k$. \hfill{$\square$}

\vspace{0.2cm}

Note that Jans showed in [J2, Theorem 4.2] that for a positive
integer $k$, fin.dim$\Lambda \leq k$ if and only if a module $N
\in$mod $\Lambda ^{op}$ of $D$-class $k$ is projective provided
$N^*$ is projective. Summarizing the results obtained above, we get
the following

\vspace{0.2cm}

{\bf Theorem 2.16} {\it Let} $k$ {\it be a non-negative integer and}
$\Lambda$ {\it a} $(k+1)$-{\it Gorenstein ring. Then the following
statements are equivalent.}

(1) id$_{\Lambda}\Lambda \leq k$.

(2) Fin.dim$\Lambda \leq k$.

(3) fin.dim$\Lambda \leq k$.

(4) {\it Each module in} mod $\Lambda ^{op}$ {\it with reduced grade
at least} $k+1$ {\it is reflexive}.

(5) {\it Each module in} mod $\Lambda ^{op}$ {\it with reduced grade
at least} $k+1$ {\it is torsionless}.

(6) {\it Each module in} mod $\Lambda ^{op}$ {\it of} $D$-{\it
class} $k$ {\it is reflexive}.

(7) {\it A module} $N \in$mod $\Lambda ^{op}$ {\it of} $D$-{\it
class} $k$ {\it is projective provided} $N^*$ {\it is projective}.

\vspace{0.2cm}

For any positive integer $n$, Zimmermann Huisgen showed in [Zi] that
there exists a finite-dimensional monomial relation algebra
$\Lambda$ with $J^4=0$ (where $J$ is the Jacobson radical of
$\Lambda$) such that fin.dim$\Lambda=n+1$ and Fin.dim$\Lambda=n+2$.
In [S] Smal$\phi$ constructed examples of algebras $\Lambda$ with
arbitrarily big difference between fin.dim$\Lambda$ and
Fin.dim$\Lambda$. But very little appears to be known when
fin.dim$\Lambda=$Fin.dim$\Lambda$. As an immediate consequence of
Theorem 2.16 we have

\vspace{0.2cm}

{\bf Corollary 2.17}
fin.dim$\Lambda=$Fin.dim$\Lambda=$id$_{\Lambda}\Lambda$ {\it for an}
{\it $\infty$-Gorenstein ring} $\Lambda$.

\vspace{0.2cm}

{\bf Corollary 2.18} {\it Let} $\Lambda$ {\it be a left and right
Artinian ring. If} $\Lambda$ {\it is} $\infty$-{\it Gorenstein,
then} fin.dim$\Lambda=$fin.dim$\Lambda ^{op}=$Fin.dim$\Lambda=
$Fin.dim$\Lambda ^{op}$=id$_{\Lambda}\Lambda=$ id$_{\Lambda
^{op}}\Lambda$.

\vspace{0.2cm}

{\it Proof.} By Corollaries 2.8 and 2.17. \hfill{$\square$}

\vspace{0.2cm}

Recall that $\Lambda$ is said to have {\it dominant dimension} at
least $k$ if each $I^{'}_i$ is flat for any $0 \leq i \leq k-1$.

The following are some important homological conjectures in
representation theory of algebras.

Nakayama Conjecture ({\bf NC})$^{\rm [Y]}$: An Artinian algebra
$\Lambda$ is self-injective if $\Lambda$ has infinite dominant
dimension.

It is easy to see that the Nakayama conjecture is a special case of
the following conjecture.

Auslander-Reiten Conjecture ({\bf ARC})$^{{\rm [AR3]}}$: An
$\infty$-Gorenstein Artinian algebra has finite left and right
self-injective dimensions.

Finitistic Dimension Conjecture ({\bf FDC})$^{\rm [B1, Zi]}$:
fin.dim$\Lambda$ is finite for an Artinian algebra $\Lambda$.

It follows from Corollary 2.18 that {\bf ARC} and {\bf FDC} are
equivalent for $\infty$-Gorenstein Artinian algebras. So we get the
following implications: {\bf FDC}$\Rightarrow${\bf ARC}$\Rightarrow$
{\bf NC}. Furthermore, we get some equivalent versions of {\bf NC}
as follows.

\vspace{0.2cm}

{\bf Corollary 2.19} {\it Let} $\Lambda$ {\it be an Artinian algebra
with infinite dominant dimension. Then the following statements are
equivalent.}

(1) $\Lambda$ {\it is self-injective}.

(2) fin.dim$\Lambda=0$.

(3) {\it The right annihilator of any proper left ideal of}
$\Lambda$ {\it is non-zero}.

(4) {\it The right annihilator of any maximal left ideal of}
$\Lambda$ {\it is non-zero}.

(5) grade$M<\infty$ {\it for any non-zero module} $M\in$mod
$\Lambda$.

(6) grade$S<\infty$ {\it for any simple module} $S\in$mod $\Lambda$.

\vspace{0.2cm}

{\it Proof.} $(1)\Leftrightarrow (2)$ follows from Corollary 2.17.
By Corollary 2.18 and [B1, Corollary 5.6 and Theorem 5.4], we have
that $(3)\Leftrightarrow (2)\Rightarrow (5)$. Both $(3)\Rightarrow
(4)$ and $(5)\Rightarrow (6)$ are trivial.

$(4)\Rightarrow (3)$ Assume that $I$ is a proper left ideal of
$\Lambda$. Then $I$ is contained in a maximal left ideal $L$ of
$\Lambda$. So the right annihilator of $L$ is contained in that of
$I$.

$(6)\Rightarrow (4)$ Assume that $S\in$mod $\Lambda$ is a simple
module. Then grade$S<\infty$ by (6). Put $K_i=$Ker$(I^{'}_i \to
I^{'}_{i+1})$ for any $i\geq 0$. Then we get an exact sequence:
$$0\to {\rm Hom}_{\Lambda}(S, K_i) \to
{\rm Hom}_{\Lambda}(S, I^{'}_i) \to {\rm Hom}_{\Lambda}(S, K_{i+1})
\to {\rm Ext}_{\Lambda}^{i+1}(S, \Lambda) \to 0$$ for any $i\geq 0$.
Since $\Lambda$ has infinite dominant dimension, $I^{'}_i$ is
projective for any $i\geq 0$.

We claim that $S^*\neq 0$. Otherwise, if $S^*=0$, then ${\rm
Hom}_{\Lambda}(S, I^{'}_i)=0$ and ${\rm Hom}_{\Lambda}(S, K_i)=0$
for any $i\geq 0$. Then by the exactness of the above sequence, we
get that ${\rm Ext}_{\Lambda}^i(S, \Lambda)=0$ for any $i\geq 0$,
which contradicts with grade$S<\infty$. The claim is proved. Now
assume that $L$ is a maximal left ideal of $\Lambda$. Then $\Lambda
/L\in$mod $\Lambda$ is a simple module. By [B1, (4.1)], we have that
the right annihilator of $L$ is isomorphic to $(\Lambda /L)^*$,
which is non-zero by the above claim. \hfill{$\square$}

\vspace{0.2cm}

The generalized Nakayama conjecture ({\bf GNC}), posed by Auslander
and Reiten in [AR1], has an equivalent version as follows: Over an
Artinian algebra $\Lambda$, grade$S<\infty$ for any simple module
$S\in$mod $\Lambda$. It is well known that {\bf
GNC}$\Rightarrow${\bf NC} (see [AR1]). Corollary 2.19 not only gives
another proof of this implication, but also shows that in order to
verify {\bf NC} it suffices to verify {\bf GNC} for Artinian
algebras with infinite dominant dimension.

\vspace{0.5cm}

\centerline{\large \bf 3. Some properties of grade of modules}

\vspace{0.2cm}

In this section, we study the properties of grade of modules over
quasi $\infty$-Gorenstein rings. We begin with the following lemma.

\vspace{0.2cm}

{\bf Lemma 3.1} ([H2, Lemma 6.2]) {\it Let} $M\in$mod $\Lambda$ {\it
and} $n$ {\it a non-negative integer. If} grade$M\geq n$ {\it and}
gradeExt$_{\Lambda}^n(M, \Lambda)\geq n+1$, {\it then} grade$M\geq
n+1$.

\vspace{0.2cm}

The following is the main result in this section.

\vspace{0.2cm}

{\bf Theorem 3.2} {\it Let} $\Lambda$ {\it be a left quasi
$\infty$-Gorenstein ring and} $M\in$mod $\Lambda$ {\it with}
grade$M$ {\it finite. Then} Ext$_{\Lambda}^i($Ext$_{\Lambda
^{op}}^i($Ext$_{\Lambda}^{{\rm grade}M}(M, \Lambda), \Lambda),
\Lambda)=0$ {\it if and only if} $i\neq$grade$M$.

\vspace{0.2cm}

{\it Proof.} Because $\Lambda$ is left quasi $\infty$-Gorenstein,
s.gradeExt$_{\Lambda ^{op}}^{i+1}(N, \Lambda)\geq i$ for any
\linebreak
$N\in$ mod $\Lambda ^{op}$ and $i\geq 1$ by Lemma 2.1.

Let $M\in$mod $\Lambda$ and
$$\cdots \to P_i \buildrel {f_i} \over \longrightarrow \cdots \buildrel {f_2}
\over \longrightarrow P_1 \buildrel {f_1} \over \longrightarrow P_0
\to M \to 0 \eqno{(1)}$$ be a projective resolution of $M$ in mod
$\Lambda$. Put $t=$grade$M$ and $X_i=$Coker$f_i^*$ for any $i\geq
1$.

If $t=0$ (that is, $M^*\neq 0$), then $M^{***} \cong M^{*} (\neq 0)$
by [H2, Lemma 1.6]. On the other hand, for any $i\geq 1$, we have
Ext$_{\Lambda ^{op}}^i(M^*, \Lambda)\cong$Ext$_{\Lambda
^{op}}^{i+2}(X_1, \Lambda)$. So gradeExt$_{\Lambda ^{op}}^i(M^*,
\Lambda)=$ \linebreak gradeExt$_{\Lambda ^{op}}^{i+2}(X_1,
\Lambda)\geq i+1$ and hence Ext$_{\Lambda}^i($Ext$_{\Lambda
^{op}}^i(M^*, \Lambda), \Lambda)=0$ for any $i\geq 1$.

Now suppose $t\geq 1$. Then from the exact sequence (1), we get the
following exact sequence:
$$0\to M^*(=0) \to P_0^* \buildrel {f_1^*} \over \longrightarrow
P_1^* \buildrel {f_2^*} \over \longrightarrow \cdots \buildrel
{f_t^*} \over \longrightarrow P_t^* \to X_t \to 0$$ in mod $\Lambda ^{op}$,
which implies
pd$_{\Lambda ^{op}}X_t \leq t$ and induces an exact sequence:
$$P_t^{**} \buildrel {f_t^{**}} \over \longrightarrow \cdots \buildrel {f_2^{**}}
\over \longrightarrow P_1^{**} \buildrel {f_1^{**}} \over
\longrightarrow P_0^{**} \to {\rm Ext}_{\Lambda}^t(X_t, \Lambda) \to
0.$$ So we have $M\cong$Ext$_{\Lambda ^{op}}^t(X_t, \Lambda)$.

By [Hu2, Lemma 2], we have an exact sequence:
$$0\to {\rm Ext}_{\Lambda}^t(M, \Lambda)\to X_t \to P_{t+1}^* \to
X_{t+1}\to 0 \eqno{(2)}$$ Put $Y_t=$Im$(X_t \to P_{t+1}^*)$. Notice
that gradeExt$_{\Lambda}^t(M, \Lambda)\geq t$ by Lemma 2.1, we then
get an exact sequence:
$$0={\rm Ext}_{\Lambda ^{op}}^{t-1}({\rm Ext}_{\Lambda}^t(M, \Lambda),
\Lambda)\to {\rm Ext}_{\Lambda ^{op}}^t(Y_t, \Lambda) \to {\rm
Ext}_{\Lambda ^{op}}^t(X_t, \Lambda) \to {\rm Ext}_{\Lambda
^{op}}^t({\rm Ext}_{\Lambda}^t(M, \Lambda), \Lambda)$$$$\to {\rm
Ext}_{\Lambda ^{op}}^{t+1}(Y_t, \Lambda) \to {\rm Ext}_{\Lambda
^{op}}^{t+1}(X_t, \Lambda)=0 \eqno{(3)}$$ and the isomorphisms:
$${\rm Ext}_{\Lambda ^{op}}^i({\rm Ext}_{\Lambda}^t(M, \Lambda),
\Lambda)\cong {\rm Ext}_{\Lambda ^{op}}^{i+1}(Y_t, \Lambda)\ {\rm
for \ any\ }i\geq t+1,$$ and $${\rm Ext}_{\Lambda ^{op}}^{i+1}(Y_t,
\Lambda) \cong {\rm Ext}_{\Lambda ^{op}}^{i+2}(X_{t+1}, \Lambda)\
{\rm for \ any\ } i\geq 0.$$ So, from the exact sequence (3), we get
an exact sequence:
$$0\to {\rm
Ext}_{\Lambda ^{op}}^{t+1}(X_{t+1}, \Lambda) \to M \buildrel {\pi
_M} \over \longrightarrow {\rm Ext}_{\Lambda ^{op}}^t({\rm
Ext}_{\Lambda}^t(M, \Lambda), \Lambda) \to {\rm Ext}_{\Lambda
^{op}}^{t+2}(X_{t+1}, \Lambda) \to 0 \eqno{(4)}$$ and an
isomorphism:
$${\rm
Ext}_{\Lambda ^{op}}^i({\rm Ext}_{\Lambda}^t(M, \Lambda),
\Lambda)\cong {\rm Ext}_{\Lambda ^{op}}^{i+2}(X_{t+1}, \Lambda) {\rm
\ for \ any}\ i\geq t+1.$$ It follows that grade${\rm Ext}_{\Lambda
^{op}}^i({\rm Ext}_{\Lambda}^t(M, \Lambda), \Lambda)=$grade${\rm
Ext}_{\Lambda ^{op}}^{i+2}(X_{t+1}, \Lambda)\geq i+1$ for any $i\geq
t+1$. Then Ext$_{\Lambda}^i($Ext$_{\Lambda
^{op}}^i($Ext$_{\Lambda}^t(M, \Lambda), \Lambda), \Lambda)=0$ for
any $i\geq t+1$. On the other hand, since grade$_{\Lambda}^t(M,
\Lambda)\geq t$, ${\rm Ext}_{\Lambda ^{op}}^i({\rm
Ext}_{\Lambda}^t(M, \Lambda), \Lambda)=0$ for any $0\leq i \leq
t-1$. So we conclude that Ext$_{\Lambda}^i($Ext$_{\Lambda
^{op}}^i($Ext$_{\Lambda}^t(M, \Lambda), \Lambda), \Lambda)=0$ for
any $i\neq t$.

Because gradeExt$_{\Lambda ^{op}}^{t+1}(X_{t+1}, \Lambda)\geq t$ and
gradeExt$_{\Lambda ^{op}}^{t+2}(X_{t+1}, \Lambda)\geq t+1$, from the
exact sequence (4) we get that Ext$_{\Lambda}^{i}($Ext$_{\Lambda
^{op}}^t($Ext$_{\Lambda}^t(M, \Lambda), \Lambda),
\Lambda)\cong$Ext$_{\Lambda}^{i}({\rm Im}\pi _M,
\Lambda)\cong$Ext$_{\Lambda}^{i}(M, \Lambda)=0$ for any $0 \leq i
\leq t-1$ (note: grade$M=t$). We claim that
Ext$_{\Lambda}^t($Ext$_{\Lambda ^{op}}^t($Ext$_{\Lambda}^t(M,
\Lambda), \Lambda), \Lambda)\neq 0$. Otherwise, we have that
gradeExt$_{\Lambda ^{op}}^t($Ext$_{\Lambda}^t(M, \Lambda),
\Lambda)\geq t+1$ by the above argument. Since
gradeExt$_{\Lambda}^t(M, \Lambda)\geq t$,
gradeExt$_{\Lambda}^t(M, \Lambda)\geq t+1$ and grade$M\geq t+1$
by Lemma 3.1, which contradicts with
grade$M=t$. The proof is finished. \hfill{$\square$}

\vspace{0.2cm}

Recall from [Bj] that an $\infty$-Gorenstein ring is called {\it
Auslander-Gorenstein} if it has finite left and right self-injective
dimensions. The following result was proved by Bj\"ork in [Bj,
Proposition 1.6] when $\Lambda$ is Auslander-Gorenstein.

\vspace{0.2cm}

{\bf Corollary 3.3} {\it Let} $\Lambda$ {\it be a left quasi
$\infty$-Gorenstein ring and} $M\in$mod $\Lambda$ {\it with}
grade$M$ {\it finite. Then} gradeExt$_{\Lambda}^{{\rm grade}M}(M,
\Lambda)=$grade$M$.

\vspace{0.2cm}

{\it Proof.} Suppose grade$M=k(<\infty)$. Since $\Lambda$ is a left
quasi $\infty$-Gorenstein ring, \linebreak gradeExt$_{\Lambda}^k(M,
\Lambda)\geq k$ by Lemma 2.1. On the other hand,
Ext$_{\Lambda}^k($Ext$_{\Lambda ^{op}}^k($Ext$_{\Lambda}^k(M,
\Lambda),\Lambda), \Lambda)\neq 0$ by Theorem 3.2. So Ext$_{\Lambda
^{op}}^k($Ext$_{\Lambda}^k(M, \Lambda), \Lambda) \neq 0$ and hence
gradeExt$_{\Lambda}^k(M, \Lambda)\leq k$. The proof is finished.
\hfill{$\square$

\vspace{0.2cm}

In viewing of the proof of Theorem 3.2, we get the
following

\vspace{0.2cm}

{\bf Corollary 3.4} {\it Let} $\Lambda$ {\it be a left quasi
$\infty$-Gorenstein ring with} id$_{\Lambda ^{op}}\Lambda
=t(<\infty)$. {\it If} $M\in$mod $\Lambda$ {\it with} grade$M=t$,
{\it then} $M \cong$Ext$_{\Lambda ^{op}}^t($Ext$_{\Lambda}^{t}(M,
\Lambda), \Lambda)$.

\vspace{0.2cm}

{\it Proof.} Consider the exact sequence (4) in the proof of Theorem
3.2. If id$_{\Lambda ^{op}}\Lambda =t$, then ${\rm Ext}_{\Lambda
^{op}}^{t+1}(X_{t+1}, \Lambda)=0={\rm Ext}_{\Lambda
^{op}}^{t+2}(X_{t+1}, \Lambda)$. So $M \cong$Ext$_{\Lambda
^{op}}^t($Ext$_{\Lambda}^{t}(M, \Lambda), \Lambda)$.
\hfill{$\square$

\vspace{0.2cm}

A left (resp. right) quasi $\infty$-Gorenstein ring is called {\it
left (resp. right) quasi Auslander-Gorenstein} if it has finite left
and right self-injective dimensions. We denote $\mathcal{G}_t= \{
M\in {\rm mod}\ \Lambda | {\rm grade}M=t\}$. The following corollary
generalizes [I2, Theorem 4].

\vspace{0.2cm}

{\bf Corollary 3.5} {\it Let} $\Lambda$ {\it be a left and right
quasi Auslander-Gorenstein ring with}
id$_{\Lambda}\Lambda$=id$_{\Lambda ^{op}}\Lambda$ \linebreak $=t$.
{\it If} $M\in$mod $\Lambda$ {\it with} grade$M=t$, {\it then} $M
\cong$Ext$_{\Lambda ^{op}}^t($Ext$_{\Lambda}^{t}(M, \Lambda),
\Lambda)$. {\it Moreover, the functors} Ext$_{\Lambda}^t(\ ,
\Lambda)$ {\it and} Ext$_{\Lambda ^{op}}^t(\ , \Lambda)$ {\it give a
duality between} $\mathcal{G}_t$ and $\mathcal{G}_t^{op}$.

\vspace{0.2cm}

{\it Proof.} By Corollaries 3.3 and 3.4. \hfill{$\square$

\vspace{0.2cm}

{\bf Example 3.6} There exist rings which are left and right quasi
Auslander-Gorenstein, but not Auslander-Gorenstein. For example, let
$\Lambda$ be the path algebra given by the quiver $2 \leftarrow 1
\rightarrow 3$. Then id$_{\Lambda}\Lambda =$id$_{\Lambda
^{op}}\Lambda =1$ and fd$_{\Lambda}I_0^{'} =$fd$_{\Lambda ^{op}}I_0
=1$. So $\Lambda$ is left and right quasi Auslander-Gorenstein, but
not Auslander-Gorenstein.

\vspace{0.2cm}

A module $M\in$mod $\Lambda$ is called {\it pure} if
grade$X=$grade$M$ for any non-zero submodule $X$ of $M$ (see [Bj]). We
use $\mathcal{C}_{\Lambda}^n$ to denote the full subcategory of mod
$\Lambda$ consisting of the modules $M$ with Hom$_{\Lambda}(M,
\bigoplus _{i=0}^n I'_i)=0$ (see [CSS]).

\vspace{0.2cm}

{\bf Lemma 3.7} {\it Let} $\Lambda$ {\it be an Auslander-Gorenstein
ring and} $k$ {\it a positive integer. Then the following statements
are equivalent for a module} $M\in$mod $\Lambda$ {\it with}
grade$M=k$.

(1) $M$ {\it is pure}.

(2)
$M\in\mathcal{C}_{\Lambda}^{k-1}\setminus\mathcal{C}_{\Lambda}^k$.

(3) Ext$_{\Lambda ^{op}}^i($Ext$_{\Lambda}^i(M, \Lambda),
\Lambda)=0$ {\it for every} $i\neq k$.

\vspace{0.2cm}

{\it Proof.} $(1)\Leftrightarrow (3)$ See [Bj, Proposition 1.9].

By [CSS, Lemma 1.1], for a non-negative integer $n$, we have that a
module $M\in$mod $\Lambda$ is in $\mathcal{C}_{\Lambda}^n$ if and
only if s.grade$M\geq n+1$. From this fact it is easy to get
$(1)\Leftrightarrow (2)$. \hfill{$\square$}

\vspace{0.2cm}

Bj\"ork raised in [Bj, p.144] a question as follows: For an
Auslander-Gorenstein ring $\Lambda$, is it true that
Ext$_{\Lambda}^{{\rm grade}M}(M, \Lambda)$ is pure for any $M\in$
mod $\Lambda$? It was answered affirmatively
by Bj\"ork and Ekstr\"om in [BjE, Proposition 2.11]. As an
application of Theorem 3.2, we give a different proof of this
result.

\vspace{0.2cm}

{\bf Proposition 3.8} {\it Let} $\Lambda$ {\it be an
Auslander-Gorenstein ring. Then} Ext$_{\Lambda}^{{\rm grade}M}(M,
\Lambda)$ {\it is pure for any} $M\in$ mod $\Lambda$.

\vspace{0.2cm}

{\it Proof.} Suppose $M\in$ mod $\Lambda$. Since $\Lambda$ has
finite self-injective dimensions, grade$M$ is finite. By Corollary
3.3, gradeExt$_{\Lambda}^{{\rm grade}M}(M, \Lambda)$=grade$M$. On
the other hand, by Theorem 3.2, we have that
Ext$_{\Lambda}^i($Ext$_{\Lambda ^{op}}^i($Ext$_{\Lambda}^{{\rm
grade}M}(M, \Lambda), \Lambda), \Lambda)=0$ if and only if
$i\neq$grade$M$. It follows from Lemma 3.7 that Ext$_{\Lambda}^{{\rm
grade}M}(M, \Lambda)$ is pure. \hfill{$\square$}

\vspace{0.2cm}

In view of the results obtained above, it is natural to ask the
following question.

\vspace{0.2cm}

{\bf Question} {\it Does Proposition 3.8 hold true for left (and
right) quasi Auslander-Gorenstein rings? That is, for a left (and
right) quasi Auslander-Gorenstein ring $\Lambda$, is it true that
${\rm Ext}_{\Lambda}^{{\rm grade}M}(M, \Lambda)$ is pure for any
$M\in{\rm mod}\ \Lambda$?}

\vspace{0.2cm}

{\it Remark.} This question is a generalized version of the
Bj\"ork's question above. By the proof of Proposition 3.8, it is
easy to see that the answer to this question is affirmative if the
implication $(3)\Rightarrow (1)$ in Lemma 3.7 also holds true for
left (and right) quasi Auslander-Gorenstein rings.

\vspace{0.2cm}

In the following, we give some further properties of grade of
modules.

Let $0\to M_1 \to M_2 \to M_3 \to 0$ be an exact sequence in mod
$\Lambda$. In general, we have grade$M_2 \geq$ min$\{$grade$M_1$,
grade$M_3 \} $. Bj\"ork showed in [Bj, Proposition 1.8] that the
equality holds true if $\Lambda$ is an Auslander-Gorenstein ring.
The following result shows that the assumption ``$\Lambda$ is
Gorenstein" is not necessary for this Bj\"ork's result.

\vspace{0.2cm}

{\bf Proposition 3.9} {\it Let} $\Lambda$ {\it be an
$\infty$-Gorenstein ring and} $0\to M_1 \to M_2 \to M_3 \to 0$ {\it
an exact sequence in} mod $\Lambda$. {\it Then} grade$M_2
=$min$\{$grade$M_1$, grade$M_3 \}$.

\vspace{0.2cm}

{\it Proof.} It suffices to prove grade$M_2 \leq$min$\{$grade$M_1$,
grade$M_3 \} $. Put $n=$min$\{$grade$M_1$, grade$M_3 \} $. Without
loss of generality, suppose $n<\infty$. We proceed in three cases.

{\it Case I.} Assume that $n=$grade$M_1 =$grade$M_3$.

Consider the following exact sequence:
$$0={\rm Ext}_{\Lambda}^{n-1}(M_1 , \Lambda)\to
{\rm Ext}_{\Lambda}^{n}(M_3 , \Lambda)\to {\rm
Ext}_{\Lambda}^{n}(M_2 , \Lambda).$$ If ${\rm Ext}_{\Lambda}^{n}(M_2
, \Lambda)=0$, then ${\rm Ext}_{\Lambda}^{n}(M_3 , \Lambda)=0$ and
grade$M_3 \geq n+1$, which is a contradiction. So ${\rm
Ext}_{\Lambda}^{n}(M_2 , \Lambda)\neq 0$ and grade$M_2 \leq n$.

{\it Case II.} Assume that $n=$grade$M_3 <$grade$M_1$.

Consider the following exact sequence:
$$0={\rm Ext}_{\Lambda}^{n-1}(M_1 , \Lambda)\to
{\rm Ext}_{\Lambda}^{n}(M_3 , \Lambda)\to {\rm
Ext}_{\Lambda}^{n}(M_2 , \Lambda) \to {\rm Ext}_{\Lambda}^{n}(M_1 ,
\Lambda)=0.$$ So ${\rm Ext}_{\Lambda}^{n}(M_2 , \Lambda) \cong {\rm
Ext}_{\Lambda}^{n}(M_3 , \Lambda)\neq 0$ and hence grade$M_2 \leq
n$.

{\it Case III.} Assume that $n=$grade$M_1 <$grade$M_3$.

Consider the following exact sequence:
$$0={\rm
Ext}_{\Lambda}^{n}(M_3 , \Lambda)\to {\rm Ext}_{\Lambda}^{n}(M_2 ,
\Lambda)\to {\rm Ext}_{\Lambda}^{n}(M_1 , \Lambda)\to {\rm
Ext}_{\Lambda}^{n+1}(M_3 , \Lambda).$$ If ${\rm
Ext}_{\Lambda}^{n}(M_2 , \Lambda)=0$, then ${\rm
Ext}_{\Lambda}^{n}(M_1 , \Lambda)$ is isomorphic to a submodule of
${\rm Ext}_{\Lambda}^{n+1}(M_3 , \Lambda)$. Since $\Lambda$ is an
$\infty$-Gorenstein ring, grade${\rm Ext}_{\Lambda}^{n}(M_1 ,
\Lambda)\geq n+1$. By Lemma 3.1, we have that grade$M_1 \geq n+1$,
which is a contradiction. So ${\rm Ext}_{\Lambda}^{n}(M_2 ,
\Lambda)\neq 0$ and grade$M_2 \leq n$. \hfill{$\square$}

\vspace{0.2cm}

For a positive integer $k$, recall again that a module $M\in$mod
$\Lambda$ is called {\it $k$-torsionfree} if r.grade$\Tr M \geq
k+1$. We use $\mathcal{T}^k({\rm mod}\ \Lambda)$ to denote the full
subcategory of mod $\Lambda$ consisting of $k$-torsionfree modules.
Recall that a full subcategory $\mathcal{X}$ of mod $\Lambda$ is
said to be {\it closed under extensions} if the middle term $B$ of
any short sequence $0 \to A \to B \to C \to 0$ is in $\mathcal{X}$
provided that the end terms $A$ and $C$ are in $\mathcal{X}$. By
Lemma 2.1 and the proof of [Hu1, Lemma 3.2], we have that $\Lambda$
is a right quasi $k$-Gorenstein ring if and only if
$\mathcal{T}^i({\rm mod}\ \Lambda)$ is closed under extensions for
any $1\leq i \leq k$. The assertion (2) in the following proposition
can be regarded as a generalization of this result.

\vspace{0.2cm}

{\bf Proposition 3.10} {\it Let} $0\to M_1 \to M_2 \to M_3 \to 0$
{\it be an exact sequence in} mod $\Lambda$. {\it For a positive
integer} $k$, {\it if} grade$C\geq k$ {\it where} $C=$Coker$(M_2^*
\to M_1^*)$, {\it then we have}

(1) {\it If} $M_2 \in \mathcal{T}^{k+1}({\rm mod}\ \Lambda)$ {\it
and} $M_3 \in \mathcal{T}^k({\rm mod}\ \Lambda)$, {\it then} $M_1
\in \mathcal{T}^{k+1}({\rm mod}\ \Lambda)$.

(2) {\it If} $M_1$, $M_3\in \mathcal{T}^k({\rm mod}\ \Lambda)$, {\it
then} $M_2\in \mathcal{T}^k({\rm mod}\ \Lambda)$.

(3) {\it If} $M_1 \in \mathcal{T}^k({\rm mod}\ \Lambda)$ {\it and}
$M_2 \in \mathcal{T}^{k-1}({\rm mod}\ \Lambda)$, {\it then} $M_3 \in
\mathcal{T}^{k-1}({\rm mod}\ \Lambda)$.

\vspace{0.2cm}

{\it Proof.} Consider the following commutative diagram with exact
columns and rows:
$$\xymatrix{& 0 & 0 & 0 &\\
0\ar[r] & M_1\ar[u]\ar[r] & M_2\ar[u]\ar[r] & M_3\ar[u]\ar[r]
& 0\\
0\ar[r] & F_{0}\ar[u]\ar[r] & F_{0}\oplus G_{0}\ar[u]\ar[r] &
G_{0}\ar[u]\ar[r] & 0\\
0\ar[r] & F_{1}\ar[u]\ar[r] & F_{1}\oplus G_{1}\ar[u]\ar[r] &
G_{1}\ar[u]\ar[r] & 0 }$$ where all $F_{i}$ and $G_{i}$ are
projective in mod $\Lambda$. Then we get the following commutative
diagram with exact columns and rows:
$$\xymatrix{& 0\ar[d] & 0\ar[d] & 0\ar[d] &\\
0\ar[r] & M_3^{*}\ar[d]\ar[r] & M_2^{*}\ar[d]\ar[r] &
M_1^{*}\ar[d] & \\
0\ar[r] & G_{0}^{*}\ar[d]\ar[r] & G_{0}^{*}\oplus
F_{0}^{*}\ar[d]\ar[r] & F_{0}^{*}\ar[d]\ar[r] & 0\\
0\ar[r] & G_{1}^{*}\ar[d]\ar[r] & G_{1}^{*}\oplus
F_{1}^{*}\ar[d]\ar[r] & F_{1}^{*}\ar[d]\ar[r] & 0\\
& \Tr M_3\ar[d] & \Tr M_2\ar[d] & \Tr M_1\ar[d] &\\
& 0 & 0 & 0 & }$$ It follows from the snake lemma that $0\to
M_3^{*}\to M_2^{*} \to M_1^{*} \to \Tr M_3\to \Tr M_2\to \Tr M_1\to
0$ is exact. Then we get two short exact sequences:
$0\to C\to \Tr M_3 \to K\to 0$ and $0\to K\to \Tr M_2 \to \Tr M_1\to
0$, where $C=$Ker$(\Tr M_3 \to \Tr M_2)$ and $K=$Im$(\Tr M_3 \to \Tr
M_2)$, which yield two long exact sequences:

$${\rm Ext}_{\Lambda ^{op}}^i(K, \Lambda)\to {\rm
Ext}_{\Lambda ^{op}}^i(\Tr M_3, \Lambda)\to {\rm Ext}_{\Lambda
^{op}}^i(C, \Lambda)\to {\rm Ext}_{\Lambda ^{op}}^{i+1}(K,
\Lambda)$$$$\to {\rm Ext}_{\Lambda ^{op}}^{i+1}(\Tr M_3, \Lambda)\to
{\rm Ext}_{\Lambda ^{op}}^{i+1}(C, \Lambda)$$ and
$${\rm Ext}_{\Lambda ^{op}}^i(\Tr M_1, \Lambda)\to {\rm
Ext}_{\Lambda ^{op}}^i(\Tr M_2, \Lambda)\to {\rm Ext}_{\Lambda
^{op}}^i(K, \Lambda)\to {\rm Ext}_{\Lambda ^{op}}^{i+1}(\Tr M_1,
\Lambda)$$$$\to {\rm Ext}_{\Lambda ^{op}}^{i+1}(\Tr M_2, \Lambda)\to
{\rm Ext}_{\Lambda ^{op}}^{i+1}(K, \Lambda)$$ for any $i\geq 0$.

If $M_2\in \mathcal{T}^{k+1}({\rm mod}\ \Lambda)$ and $M_3\in
\mathcal{T}^k({\rm mod}\ \Lambda)$, then r.grade$\Tr M_2\geq k+2$
and r.grade$\Tr M_3$
\linebreak
$\geq k+1$. Because grade$C\geq k$
by assumption, from the above two long exact sequences we get that
r.grade$K\geq k+1$ and ${\rm Ext}_{\Lambda ^{op}}^i(\Tr M_1,
\Lambda)=0$ for any $2\leq i \leq k+1$. But $M_1$ is clearly
torsionless, so ${\rm Ext}_{\Lambda ^{op}}^1(\Tr M_1, \Lambda)=0$
and r.grade$\Tr M_1\geq k+2$, which implies that $M_1$ is
$(k+1)$-torsionfree. This finishes the proof of (1). Similarly, we
get (2) and (3). \hfill{$\square$}

\vspace{0.2cm}

We end this section by giving some examples of rings satisfying the
grade condition ``grade$C\geq k$" in Proposition 3.10.

\vspace{0.2cm}

{\bf Example 3.11} (1) From the proof of [Hu1, Theorem 2.3], we know
that if fd$_{\Lambda ^{op}}\bigoplus _{i=0}^{k-1}I_i\leq k$, then
the grade condition in Proposition 3.10 is satisfied. In particular,
if $\Lambda$ is a right quasi $k$-Gorenstein ring, then this grade
condition is also satisfied.

(2) By [I1, Proposition 1], we have that
id$_{\Lambda}\Lambda=$sup$\{$fd$_{\Lambda ^{op}}I|I$ is an injective
right $\Lambda$-module$\}$. Then by (1), the grade condition in
Proposition 3.10 is satisfied if id$_{\Lambda}\Lambda \leq k$. Thus,
if id$_{\Lambda}\Lambda =$id$_{\Lambda ^{op}}\Lambda \leq k$, then
by Proposition 3.10, $\mathcal{T}^k({\rm mod}\ \Lambda)$ is a
resolving subcategory of mod $\Lambda$ in the sense of Auslander and
Reiten [AR2].

\vspace{0.5cm}

\centerline{\large \bf 4. Torsionless and reflexive modules}

\vspace{0.2cm}

Let $A\in$mod $\Lambda$ be a torsionless module. Then $A$ can be
embedded into a finitely generated free $\Lambda$-module $G$. We use
$\mathcal{E}_A$ to denote the subcategory of mod $\Lambda$
consisting of the non-zero modules $C$ such that there exists an
exact sequence $0 \to A \to G \to C \to 0$ in mod $\Lambda$ with $G$
free.

\vspace{0.2cm}

{\bf Proposition 4.1} {\it If} $\Lambda$ {\it is a right quasi}
$k$-{\it Gorenstein ring, then, for any} $t$-{\it torsionfree
module} $A\in$ mod $\Lambda$ {\it (where} $1\leq t \leq k${\it )}
{\it and} $C\in \mathcal{E}_A$, {\it there exists an exact sequence}
$0 \to F \to T \to C \to 0$ {\it in} mod $\Lambda$ {\it with} $F$
{\it free and} $T$ $(t-1)$-{\it torsionfree}.

\vspace{0.2cm}

{\it Proof.} Because $\Lambda$ is a right quasi $k$-Gorenstein ring,
we have that for any $1\leq t \leq k$, a module in mod $\Lambda$ is
$t$-torsionfree if and only if it is $t$-syzygy by [AR4, Proposition
1.6 and Theorem 1.7]. So for a $t$-torsionfree module $A\in$mod
$\Lambda$, there exists an exact sequence $0\to A \to F \to K \to 0$
in mod $\Lambda$ with $F$ free and $K$ $(t-1)$-torsionfree. Let
$C\in\mathcal{E}_A$. Then there exists an exact sequence $0\to A \to
G \to C \to 0$ in mod $\Lambda$ with $G$ free. Consider the
following push-out diagram:
$$\xymatrix{& 0 \ar[d] & 0 \ar[d] & & \\
0\ar[r] & A\ar[r] \ar[d] & G \ar[r] \ar[d] & C \ar[r]
\ar@{=}[d] & 0\\
0\ar[r] & F\ar[r] \ar[d] & T \ar[r] \ar[d] & C \ar[r] & 0\\
& K \ar@{=}[r] \ar[d] & K \ar[d] & &\\
& 0 & 0 & &}
$$ By Lemma 2.1 and the proof of [Hu1, Lemma
3.2], $\mathcal{T}^t({\rm mod}\ \Lambda)$ is closed
under extensions for any $1\leq t \leq k$. So from the exactness of
the middle column in the above diagram, we get that $T$ is $(t-1)$-torsionfree.
Hence the middle row in the above diagram is as desired. \hfill{$\square$}

\vspace{0.2cm}

Recall from [CSS] that a module $M\in$mod $\Lambda$ is said to be
{\it pseudo-null} if $M\in \mathcal{C}_{\Lambda}^1$ (that is,
Hom$_{\Lambda}(M, I'_0 \bigoplus I'_1)=0$).

\vspace{0.2cm}

{\bf Proposition 4.2} {\it Let} $\Lambda$ {\it be a 2-Gorenstein
ring, and let} $A\in$mod $\Lambda$ {\it be a torsionless module and} $0 \to A
\to G \to C \to 0$ {\it an exact sequence in} mod $\Lambda$ {\it
with} $G$ {\it free and} $C\in \mathcal{E}_A$. {\it Then we have}

(1) Coker$\sigma _A$ {\it is pseudo-null}.

(2) $A^{**}$ {\it is isomorphic to a submodule of} $G$.

(3) {\it If} $C$ {\it has no non-zero pseudo-null submodule, then}
$A$ {\it is reflexive}.

\vspace{0.2cm}

{\it Proof.} (1) By Lemma 2.12, Coker$\sigma _A \cong$Ext$_{\Lambda
^{op}}^2(\Tr A, \Lambda)$. Since $\Lambda$ is 2-Gorenstein,
\linebreak
Hom$_{\Lambda}($Ext$_{\Lambda ^{op}}^2(\Tr A, \Lambda),
I_0^{'} \bigoplus I_1^{'})=0$ by [IS, Proposition 3].

(2) From the exact sequence $0 \to A \to G \to C \to 0$ we get an
exact sequence:
$$0\to C^* \to G^* \to A^* \to {\rm Ext}_{\Lambda}^1(C,
\Lambda)\to 0.$$ Put $K=$Im$(G^* \to A^*)$. Then we get two exact
sequences $0\to K^* \to G^{**}(\cong G)$ and $[{\rm
Ext}_{\Lambda}^1(C, \Lambda)]^* \to A^{**} \to K^*$. Because
$\Lambda$ is 2-Gorenstein, s.gradeExt$_{\Lambda}^1(C, \Lambda)\geq
1$ and $[{\rm Ext}_{\Lambda}^1(C, \Lambda)]^* =0$. Then the
assertion follows.

(3) Note that $A^{**}/A\cong$Coker$\sigma _A$ is pseudo-null by
(1) and that $A^{**}/A$ is isomorphic to a submodule of $C$ by (2).
So $A^{**}/A =0$ by assumption and hence $A\cong A^{**}$ and
$A$ is reflexive. \hfill{$\square$}

\vspace{0.2cm}

We now give a criterion for judging when a torsionless module is
reflexive.

\vspace{0.2cm}

{\bf Theorem 4.3} {\it Let} $\Lambda$ {\it be a 2-Gorenstein ring
and} $A \in$mod $\Lambda$ {\it a torsionless module}. {\it Then the
following statements are equivalent.}

(1) $A$ {\it is reflexive.}

(2) $C$ {\it has no non-zero pseudo-null submodule for any}
$C\in\mathcal{E}_A$.

(3) $C$ {\it has no non-zero pseudo-null submodule for some}
$C\in\mathcal{E}_A$.

\vspace{0.2cm}

{\it Proof.} $(2)\Rightarrow (3)$ is trivial, and $(3)\Rightarrow
(1)$ follows from Proposition 4.2(3).

$(1)\Rightarrow (2)$ Assume that $A \in$mod $\Lambda$ is reflexive.
Then, by Proposition 4.1, for any $C\in\mathcal{E}_A$, there exists
an exact sequence $0 \to F \to T \to C \to 0$ in mod $\Lambda$ with
$F$ free and $T$ torsionless. By [M, Corollary 1.3], it is easy to
see that $C$ can be embedded into a finite direct sum of $I'_0
\bigoplus I'_1$ and so $C$ has no non-zero pseudo-null submodule.
\hfill{$\square$}

\vspace{0.2cm}

Let $\Lambda$ be an Auslander-Gorenstein ring and $I$ a non-zero
proper left ideal of $\Lambda$. Then, by Lemma 3.7, we have that
$\Lambda /I$ has no non-zero pseudo-null submodule if and only if
$\Lambda /I$ is pure of grade 1. It was showed in [CSS, Lemma 3.12]
that if $\Lambda$ is an Auslander-regular ring (that is, $\Lambda$ is an
$\infty$-Gorenstein ring with finite global dimension) without
non-zero zero divisors, then $I$ is reflexive if and only if
$\Lambda /I$ is pure of grade 1. The following corollary generalizes
this result, which is an immediate consequence of Theorem 4.3.

\vspace{0.2cm}

{\bf Corollary 4.4} {\it Let} $\Lambda$ {\it be a 2-Gorenstein ring.
Then a non-zero proper left ideal} $I$ {\it of} $\Lambda$ {\it is
reflexive if and only if} $\Lambda /I$ {\it has no non-zero
pseudo-null submodule}.

\vspace{0.2cm}

Recall from [FuW] that a right $\Lambda$-module $M$ is said to {\it
have an injective resolution with a redundant image} from a positive
integer $n$ if the $n$-th cosyzygy $\Omega _n$ has a decomposition
$\Omega _n=\bigoplus _{i\in I}A_i$ such that each $A_i$ is a direct
summand of a cosyzygy $\Omega _{\alpha _i}$ for some $\alpha _i\neq
n$. It was showed in [FuW, Theorem 1] that if $\Lambda _{\Lambda}$
has an injective resolution with a redundant image from $n$, then
$\bigoplus _{i=0}^n I'_i$ is an injective cogenerator for the
category of left $\Lambda$-modules.

Assume that $\Lambda$ is a 2-Gorenstein ring such that $\Lambda
_{\Lambda}$ has an injective resolution with a redundant image from
1. Then $I'_0 \bigoplus I'_1$ is an injective cogenerator for the
category of left $\Lambda$-modules. So every module in mod $\Lambda$
has no non-zero pseudo-null submodule and hence each torsionless
module in mod $\Lambda$ is reflexive by Theorem 4.3. It then follows from
[J2, Theorem 5.1] that id$_{\Lambda ^{op}}\Lambda \leq 1$. Therefore we
have established the following result.

\vspace{0.2cm}

{\bf Corollary 4.5} {\it Let} $\Lambda$ {\it be a 2-Gorenstein ring.
If} $\Lambda _{\Lambda}$ {\it has an injective resolution with a
redundant image from 1, then} id$_{\Lambda ^{op}}\Lambda \leq 1$.

\vspace{0.2cm}

Ramras raised in [G, p.380] an open question: When is each reflexive
module in mod $\Lambda$ projective? A generalized version of this
question is: For a positive integer $k$, when is each
$k$-torsionfree module in mod $\Lambda$ projective? In the
following, we will deal with these two questions and give some partial
answers to them.

\vspace{0.2cm}

{\bf Proposition 4.6} {\it For any positive integer} $k$, {\it the
following statements are equivalent.}

(1) {\it Each} $k$-{\it torsionfree module in} mod $\Lambda$ {\it is
projective.}

(2) {\it Each module in} mod $\Lambda ^{op}$ {\it with reduced grade
at least} $k+1$ {\it is projective.}

\vspace{0.2cm}

{\it Proof.} $(1)\Rightarrow (2)$ Let $N\in$mod $\Lambda ^{op}$ with
r.grade$N\geq k+1$ and $Q_1 \buildrel {g} \over \longrightarrow Q_0
\to N \to 0$ be a projective presentation of $N$ in mod $\Lambda
^{op}$. Then we get an exact sequence:
$$0 \to N^* \to Q^*_0 \buildrel {g^*} \over
\longrightarrow Q^*_1 \to \Tr N \to 0.$$
It is easy to see that $N\cong$Coker$g^{**}=\Tr\Tr N$. So $\Tr
N\in$mod $\Lambda$ is $k$-torsionfree and hence $\Tr N$ is
projective by (1). It follows from Proposition 2.13 that $N$ is
projective.

$(2)\Rightarrow (1)$ Let $M \in$mod $\Lambda$ be a $k$-torsionfree
module. Then r.grade$\Tr M\geq k+1$ and so $\Tr M$ is projective
by (2). It follows from Proposition 2.13 that $M$ is projective.
\hfill{$\square$}

\vspace{0.2cm}

It is well known that $\Lambda$ is a hereditary ring (that is,
gl.dim$\Lambda\leq 1$) if and only if each torsionless module in mod
$\Lambda$ (or mod $\Lambda ^{op}$) is projective. By Proposition
4.6, we give a new characterization of hereditary rings as follows.

\vspace{0.2cm}

{\bf Corollary 4.7} $\Lambda$ {\it is hereditary if and only if each
module in} mod $\Lambda$ ({\it or} mod $\Lambda ^{op}$) {\it with
reduced grade at least 2 is projective}.

\vspace{0.2cm}

{\bf Theorem 4.8} {\it Let} $k$ {\it be a positive integer or
infinite. Then the following statements are equivalent.}

(1) {\it Each} $k$-{\it torsionfree module in} mod $\Lambda$ {\it is
projective.}

(2) {\it Each module in} mod $\Lambda ^{op}$ {\it with reduced grade
at least} $k+1$ {\it is projective.}

\vspace{0.2cm}

{\it Proof.} When $k$ is a positive integer, it has been proved in
Proposition 4.6. When $k$ is infinite, the proof is similar to that
of Proposition 4.6, so we omit it. \hfill{$\square$}

\vspace{0.2cm}

{\bf Corollary 4.9} {\it If} gl.dim$\Lambda\leq k$, {\it then each}
$k$-{\it torsionfree module in} mod $\Lambda$ {\it is projective.}

\vspace{0.2cm}

{\it Proof.} Let $N\in$mod $\Lambda ^{op}$ with r.grade$N\geq k+1$.
Because gl.dim$\Lambda\leq k$, pd$_{\Lambda ^{op}}N\leq k$. Then by
Lemma 2.10, $N$ is projective. Thus the assertion follows from
Theorem 4.8. \hfill{$\square$}

\vspace{0.2cm}

In particular, if putting $k=2$ in Corollary 4.9, then we get the
following

\vspace{0.2cm}

{\bf Corollary 4.10} {\it If} gl.dim$\Lambda\leq 2$, {\it then each
reflexive module in} mod $\Lambda$ {\it is projective.}

\vspace{1cm}

{\bf Acknowledgements} The research was
partially supported by NSFC
(10771095, 10871088), SRFDP (200802840003), Cultivation
Fund of the Key Scientific and Technical Innovation Project, Ministry of Education of
China (No.708044) and NSF of Jiangsu Province of China.
The authors thank the referee for the
valuable suggestions.

\end{document}